\documentclass[a4paper,10pt]{article}

\usepackage[T1]{fontenc}
\usepackage[brazilian, english]{babel}
\usepackage{enumerate}
\usepackage{amsmath,amssymb,amsthm}
\usepackage{graphicx}
\usepackage{bbm}
\usepackage{dsfont}
\usepackage{cite}      
\usepackage{mathtools, nccmath, relsize}

\usepackage[hidelinks,bookmarksnumbered]{hyperref}
\hypersetup{pdfstartview=FitH}

\usepackage{tikz}
\usetikzlibrary{arrows.meta}
\usetikzlibrary{patterns}
\usetikzlibrary{shapes}
\usetikzlibrary{calc}
\usetikzlibrary{backgrounds}
\usetikzlibrary{decorations.pathreplacing}
\usetikzlibrary{decorations.pathmorphing}

  \newcounter{constant}

\def\arraypar#1{\parbox[c]{\textwidth - 2cm}{\centering #1}}

\usepackage{environ}
\NewEnviron{display}{
  \begin{equation}\begin{array}{c}
    \arraypar{\BODY}
  \end{array}\end{equation}
}

\newcommand{\Z}{{\mathbb{Z}}}

\newcommand{\N}{{\mathbb{N}}}
\newcommand{\LL}{{\mathbb{L}}}
\newcommand{\ind}{\mathbbm{1}}

\newcommand{\G}{\mathcal{G}}

\newcommand{\T}{\mathcal{T}}
\newcommand{\R}{\mathcal{R}}

\newcommand{\I}{\mathcal{I}}
\newcommand{\ta}{t^{\ast}}

\newcommand{\B}{\mathcal{B}}

\newcommand{\E}{\mathbb{E}}

\theoremstyle{plain}
\newtheorem{teo}{Theorem}[section]
\newtheorem{theorem}[teo]{Theorem}

\theoremstyle{definition}
\newtheorem{remark}{Remark}[section]

\title{Long-range contact process and percolation on a random lattice}

\author{Pablo A. Gomes\footnote{Departamento de Estat\'\i stica, Universidade de S\~ao Paulo, Rua do Mat\~ao 1010 CEP 05508-090 S\~ao Paulo-SP, Brazil. \url{pagomes@usp.br}} \and Bernardo N. B. de Lima\footnote{Departamento de Matem{\'a}tica, Universidade Federal de Minas Gerais, Av. Ant\^onio
Carlos 6627 C.P. 702 CEP 30123-970 Belo Horizonte-MG, Brazil. \url{bnblima@mat.ufmg.br}}}

\date{}

\begin{document}
\maketitle

\numberwithin{equation}{section}

\begin{abstract}
We study the phase transition phenomena for long-range oriented percolation and contact process. We study a contact process in which the range of each vertex are independent, updated dynamically  and given by some distribution $N$. We also study an analogous oriented percolation model on the hyper-cubic lattice, here there is a special direction where long-range oriented bonds are allowed; the range of all vertices are given by an i.i.d. sequence of random variables with common distribution $N$. For both models, we prove some results about the existence of a phase transition in terms of the distribution $N$.
\end{abstract}

\bigskip

\noindent{\footnotesize Keywords: contact process; long-range percolation; anisotropic percolation \\
MSC numbers:  60K35, 82B43}

\section{Introduction}

\subsection{Models and Results}

\noindent Given an integer $d \geq 1$, let $\G^d = (\Z^d, \E)$ be the complete graph whose vertex set is $\Z^d$, that is, $\E = \{\langle x,y \rangle : x \neq y\in \Z^d \}$. On this graph, we define a contact process model with two parameters, a positive real number $\lambda$, the rate of infection, and a random variable $N$ taking values on $\Z_{+}$, the range variable. Define the following sequences of random variables:
\begin{itemize}
    \item $\{R_{x,n} : x \in \Z^d, n \in \N \}$, i.i.d. with distribution $\exp(1)$;
    \item  $\{I_{e,n} : e \in \E, n \in \N\}$, i.i.d. with distribution $\exp(\lambda)$;
    \item $\{T_{x,n} : x \in \Z^d, n \in \N \}$, i.i.d. with distribution $\exp(1)$;
    \item $\{N_{x,n} : x \in \Z^d, n \in \N \}$, i.i.d. with distribution equal to that of $N$.  
\end{itemize}

We also assume that all these sequences are independent of each other. For each $x \in \Z^d$, the Poisson Process $\R_x = \{ \sum_{k=1}^n R_{x,k} : n \in \N \}$ represents the {\em recovering times} of vertex $x$; the Poisson Process $\T_x = \{ \sum_{k=1}^n T_{x,k} : n \in \N \}$ determines the {\em change times} at which the infection range of vertex $x$ is updated. 

Given $x \in \Z^d$, define
\begin{equation}\label{eq: S}
    S_{x,0} = 0 \quad  \text{and} \quad S_{x,n} = \sum_{k=1}^n T_{x,k}, ~n \geq 1
\end{equation}  
 and 
\begin{equation}\label{eq: R}
     r_x(t) = N_{x,n}, ~\forall t \in [S_{x,n-1}, S_{x,n} ), \quad n\geq1.
\end{equation}
The random variable $r_x(t)$ is the {\em range of $x$ at time $t$}.

For each bond $e \in \E$, the Poisson Process  $\I_e = \{ \sum_{k=1}^n I_{e,k} : n \in \N \}$  determines the possible {\em infection times}. At some time $t\in\I_e$, an infected end-vertex of $e$, say $x$, infects the other end-vertex, say $y$, if $|x-y|\leq r_x(t)$ and $y$ is healthy at time $t$, where $|\cdot|$ is the $\ell_1$ norm in $\Z^d$.

For each $t \geq 0$, the function $\zeta_t: \Z^d \longrightarrow \{0,1\}$ to be defined below, will denote if the vertex $x \in \Z^d$ is \emph{infected} at time $t$ if $\zeta_t(x) = 1$, and \emph{healthy} if $\zeta_t(x) = 0$. The set of infected vertices at time $t$ is denoted by $\zeta_t = \{ x \in \Z^d : \zeta_t(x) = 1\}$. 

Let $o$ be the origin of $\Z^d$. At time $t=0$, we consider the initial condition $\zeta_0$, where $\zeta_0(o)=1$ and $\zeta_0(x)=0$ for all $x\in\Z^d\setminus\{o\}$, that is, the model starts with an unique infected vertex, $o$. This choice is natural since $\Z^d$ is a transitive graph and due the additivity property (see page 32 of~\cite{Li2}).

Given any $t > 0$, the function $\zeta_t : \Z^d \longrightarrow \{0,1\}$ is defined such that $\zeta_t(y) = 1$, if and only if, for some $m \in \N$, there exist times $0 = t_0 < t_1 < \cdots < t_m \leq t$ and vertices $o = x_0, \dots, x_m = y$, such that for all  $k\in\{0, \dots, m-1\}$, setting $e_k = \langle x_k, x_{k+1} \rangle,$ it holds that:
\begin{itemize}
    \item $t_{k+1} \in \I_{e_k}$;
    \item $|e_k| \leq r_{x_k}(t_{k+1})$;
    \item $[t_k, t_{k+1}] \cap \R_{x_k} = \emptyset$;
    \item $[t_m, t] \cap \R_y = \emptyset$.
\end{itemize}

In this case, we say that there is an \emph{infection path} between $(o,0)$ and $(y,t)$. Informally, the first item above says that the infection spreads through the vertices $x_0, \dots, x_m$ at times $t_1< \cdots < t_m$; the second one says that along this infection path the spreading is restricted to the range of infection in each vertex; and the last two items say that there is no recovering times along the infection path.

Once defined the status of each vertex at any time, the construction of our processes is concluded. We refer to this stochastic processes, $(\zeta_t)_{t \geq 0}$, as Contact Processes with Dynamical Range (CPDR for short). 

Let $P$ be a probability measure under which this process is defined; although $P$ depends on $\lambda$ and $N$, we omit this from the notation. Below, we state our main results, they characterize the survival of the infection in the CPDR on $\G^d$ in terms of the distribution of $N$. 

\smallskip
\begin{theorem}\label{cfinito2}
If $E[N^d] < \infty$, then there exists $\lambda_0$ small enough such that, $\forall\ 0<\lambda<\lambda_0$, it holds that:\begin{equation}
P( \zeta_t \neq \emptyset, ~\forall t \geq 0) = 0.
\end{equation}
\end{theorem}

\smallskip
\begin{theorem}\label{cinf2}
If $\limsup_{n \to \infty} n P(N^d \geq n) > 0$, then 
\begin{equation}
P( \zeta_t \neq \emptyset, ~\forall t \geq 0) > 0, \quad \forall \lambda>0.
\end{equation}
\end{theorem}

\smallskip
It is natural to think how will behave CPDR if we consider the rates of recovering or change times different from one. Without loss of generality, we can take the rate of recovering as one (like in ordinary Contact Process) re-scaling the time variable. Concerning the rate of change times, if it is not one this will not change Theorem~\ref{cinf2}, while the statement of Theorem~\ref{cfinito2} keeps the same but the constant $\lambda_0$ can change. Note that when the rate of the change times is equals to $0$, the model is static. Unfortunately, our arguments does not extends to this case, and we believe that the hypothesis of the Theorem~\ref{cinf2} is not sufficient to guarantee positive probability of survival for all infection rate $\lambda$.

Analogous results to the theorems above can also be settled on the context of a long-range oriented percolation, where the range of bonds starting from some vertex is bounded by some random variable. More precisely, we will study the following long-range percolation model.

For any $d\geq 2$, let $(\vec{e}_i)_{i=1}^d$ be the canonical base of $\Z^d$ and $G=(\Z^d,\E_v\cup \E_h)$ be the graph where $\E_h=\cup_{n=1}^\infty \{(x,y)\in\Z^d\times\Z^d:y-x=n.\vec{e}_1\}$ and $\E_v=\{(x,y)\in\Z^d\times\Z^d:y-x=\vec{e}_i,\ i\in\{2,\dots,d\}\}$. That is, along the lines parallel to the first coordinate axis we have long-range oriented bonds and along the other directions the graph has only nearest-neighbour oriented bonds. On a random subgraph of $G$, we will define a percolation process as follows. The percolation process has three parameters $p,q\in[0,1]$ and $N$ a non-negative integer random variable. 

Given a sequence ${\bf N}=(N_x)_{x\in\Z^d}$ of i.i.d. random variables with common distribution $N$, we define the oriented random subgraph of $G$,
\begin{equation}
    G_{\bf N}:=(\Z^d,\E_v\cup (\cup_{x\in\Z^d} \{(x,x+n\vec{e}_1)\in\Z^d\times\Z^d:n\leq N_x\})).
\end{equation} 
In words, $G_{\bf N}$ is obtained from  $G$ by deleting all bonds from $x$ whose length is bigger than $N_x$. On the graph $G_{\bf N}$, we consider an independent bond percolation processes where bonds in $\E_h$ and $\E_v$ are {\em open} with probabilities $p$ and $q$, respectively (and {\em closed} with complementary probability). From now on, we call this model as Anisotropic Percolation with Random Range (APRR for short). 

Let $P$ be the underlying probability measure defined by this processes (it will be clear in each context if we are dealing with the CPDR model or the APRR model, so there is no problem in denoting the underlying measure for both models by $P$) and, as usual in percolation, we use the notation $(x\rightarrow y)$ to denote the event where there is an oriented open path from the vertex $x$ to the vertex $y$ and $(x\rightarrow \infty)$ to denote the event where the vertex $x$ is connected by oriented open paths to infinitely many vertices of $\Z^d$.

We define the percolation function and the critical curve as $\theta(p,q)=P(o\rightarrow \infty)$ and $q_c(p)=\sup\{q:\theta(p,q)=0\}$, respectively (recall that the probability measure $P$ depends on $p,q$ and $N$, thus $\theta$ and $q_c$ also depend on $N$, but for simplicity we drop these parameters in our notation).

As in the CPDR, our goal is to study how the distribution $N$ influences the behavior of the critical curve $q_c$. In this direction we prove the next two theorems.

\smallskip
\begin{theorem}\label{perc1} If $EN<\infty$, then $q_c(p)>0$ for all $p<1$.
\end{theorem}

\smallskip
In the case $EN=\infty$, we have the following partial result:

\smallskip
\begin{theorem} \label{perc2} If $\limsup_{n\rightarrow\infty} nP(N\geq n)>0$, then $q_c(p)=0$ for all $p>0$.
\end{theorem}

\smallskip
Observe that Theorem~\ref{perc2} brings up the question if under the hypothesis
$\limsup_{n\rightarrow\infty} n P(N \geq n)>0$, percolation could occur in the one-dimensional case ($q=0$). The next result is a particular case of Theorem 1 in~\cite{GGJR} and shows that even though with $q_c(p)=0$ for all $p>0$, the one-dimensional model may have zero probability
of percolation even if $p=1$.

\smallskip
\begin{theorem}[\cite{GGJR}]
\label{perc3} 
Consider the case $p=1$ and $q=0$. Let $\beta > 0$ be a constant such that the distribution of $N$ satisfies 
\begin{equation}
    P(N\geq n)=1-e^{-{\beta}/{n}}, ~ n \geq 0.
\end{equation} Then $\theta(1,0)=0$, if $\beta \leq 1$ and  $\theta(1,0)>0$, if $\beta > 1$.
\end{theorem}

\smallskip
That is, for the one-dimensional model, there is a phase transition in the parameter $\beta$. A natural question that came up from this theorem is: 
given any $0 < p < 1$, does this one-dimensional model percolate for some large $\beta$? The next result gives an affirmative answer.

\begin{theorem}\label{perc4}
Consider the case $0 < p < 1 $ and $q=0$. If $N$ has distribution as in Theorem~\ref{perc3}, then $\theta(p,0) > 0$, if $\beta > p^{-1}$.
\end{theorem}

Therefore, there is a phase transition in the parameter $\beta$ and for each $p \in [0,1]$, we have $\beta_c(p) \in [1,p^{-1}]$.

\subsection{Related Works}

\noindent The study of long-range models has its origins in the mathematical-physics literature, indeed before percolation, long-range Ising models were studied by Dyson~\cite{D1,D2} and Frolich-Spencer~\cite{FS}. Long-range percolation was first studied in one-dimensional models, where each bond $e$ with length $n$ is open, independently of each other, with probability $p_n\sim\beta n^{-s}$. Schulman~\cite{Sc} showed that there is no percolation if $s>2$. The remarkable affirmative answer was given by Newman-Schulman~\cite{NS} where it was proved that if $s<2$ there is oriented percolation; moreover, there is also non-oriented percolation in the critical case $s=2$, if $\beta$ is large enough and $p_1$ is close to one. This last result for $s=2$ was improved by Aizenman-Newman~\cite{AN}, where it was showed that there is a critical $\beta$, in the sense that there is no percolation if $\beta\leq 1$ and non-oriented percolation occurs if $\beta>1$ and $p_1$ is close to one. Afterwards, the oriented case for $s=2$ was solved by Marchetti-Sidoravicius-Vares~\cite{MSV}, also proving that $\beta_c=1$.

For long-range percolation in dimensions $d\geq 2$, an issue similar to the one treated here is the so-called {\em truncation question}: given a long-range percolation model that percolates, is there some large integer $K$ such that percolation still occurs if we delete all bonds whose lengths are bigger than $K$? The first work to tackle this question was Meester-Steif~\cite{MS} considering the case $p_n\sim e^{-cn}$. The case where $\sum p_n=\infty$ was first studied by Sidoravicius-Surgailis-Vares~\cite{SSV} for the case $p_n\sim {1}/{(n\log n)}$; Friedli-de Lima~\cite{FL} gave an affirmative answer for the case $d\geq 3$ without any additional hypothesis about $p_n$ (later, generalized for the oriented percolation and contact processes in \cite{ELV}), the case $d=2$ is still an open problem today and some partial answers had been given, see for example \cite{CL} and the references therein. One interesting negative answer was given by Biskup-Crawford-Chayes~\cite{BCC} for a long-range Potts model with $q=3$.

Long-range percolation models have been studied and showed a fruitful tool as a model for social networks, in particular the study of the graph distance on the long-range percolation cluster, the so-called {\em chemical distance}. See for example~\cite{BB,Bi1,Bi2,DS} and the references therein.

The contact process was introduced by Harris in 1974 \cite{Ha}, as a model for spreading an infection and it has became one of the most studied particle systems model since then, see the books \cite{Li1} and \cite{Li2} for an introduction on this subject. Long-range contact process is one of the several variations of this model and it had already appeared in the literature, probably its roots go back to Spitzer~\cite{Sp} and a phase transition theorem was proven in Bramson-Gray~\cite{BG}.  Long-range contact processes variations are also present in the physical literature, see \cite{GHL} and the references therein. 

More related to our work, we also would like to mention the paper of Can~\cite{Can}, in which it is studied a long-range contact processes on a percolation cluster in the one-dimensional complete graph, where each bond $\langle i,j\rangle,\ i,j\in\Z$ is open with probability $|i-j|^{-s}, s>1$.The main goal in~\cite{Can} is to give conditions on the tail of the probability of the long-range connections, such that the model still displays a positive critical infection rate, which is similar to Theorem~\ref{cfinito2}. Contact processes on a dynamical environment were the content of the recent work of Linker-Remenik \cite{LR}, it deals with a contact processes running on a one-dimensional lattice undergoing dynamical percolation, that is, the set of bonds where the contact processes is defined changes along the time, another paper concerning the same model and extending some result of~\cite{LR} for higher dimensions is~\cite{HUVV}; as in our work, the content of the these results concern about how the dynamics of the model can affect the behavior of the critical infection rate, analyzing if it is trivial or not. Other recent work concerning contact process in a dynamical environment is~\cite{SS}. Our proposal here is to merge the long-range contact process with a random and dynamical environment, in the sense of the range of infection of each individual changes along the time in a random way.

\subsection{Outline of the paper} 

\noindent Section~\ref{s2} is dedicated to the CPDR, we will prove Theorems~\ref{cfinito2} and \ref{cinf2}, in Subsections~\ref{sub21} and \ref{sub22}, respectively. 

In Section~\ref{s3}, we consider the APRR model; Theorem~\ref{perc1} is proven in Subsection~\ref{sub31}. The proof of Theorem~\ref{perc2} is similar to the one of Theorem ~\ref{cinf2}, we will point out only the differences between these proofs in Subsection~\ref{sub32}. In Subsection~\ref{sub33}, we prove the Theorem~\ref{perc4}.

Throughout this text we adopt the convention that if $x\in\Z^d$, then $x_i$ denote its $i$-th coordinate.

\section{Contact processes with dynamical range}\label{s2}

\noindent In this section, we deal with the CPDR with infection rate $\lambda$ and range distribution $N$. The demonstration of Theorem~\ref{cfinito2} is based on a standard comparison with a subcritical branching process. In the proof of Theorem~\ref{cinf2}, we will perform a block renormalization argument whose goal is to show that the CPDR dominates a supercritical, anisotropic, oriented and independent percolation model.

\subsection{Proof of Theorem \ref{cfinito2}} \label{sub21}

\noindent The idea of this proof is to define special subsets of $\Z^d\times \mathbb{R}_+$ called \emph{atoms}. These atoms will be a covering for all infection paths starting from $o$ and we will show that for $\lambda$ small enough the creation of new atoms is dominated by a subcritical branching process. This type of argument is standard to exhibit the occurrence of a subcritical phase in the contact process as well as in percolation; however, as the range variables are unbounded, a more detailed analysis is needed.

For each $x \in \Z^d$ and each $s > 0$, we define the stopping times
\begin{align}\label{eq: alphabeta}
    \alpha(x,s) &= \inf \{ t > s : t \in \R_x\}  \quad \text{and} \nonumber \\
    \beta(x,s) &= \inf \{ t > \alpha(x,s) : t \in \T_x\}.
\end{align}

The structure of atoms will be defined inductively  and the construction will guarantee that the
atoms are disjoint. Each atom is an ordered pair $A = (x,I)$, where $x \in \Z^d$ and $I \subset [0, \infty)$ is an interval. We set $s_0 = 0, a_0=0$, $x_0 = o$, $b_0 = \beta(x_0, a_0)$, $\Gamma_0 = [a_0, b_0)$ and let $A_0 = (x_0,\Gamma_0)$ be the \emph{root atom}. Observe that every part of the infection path in $\{o\}\times\Z_+$ until the first recovering time at $o$ is covered by the root atom.

Suppose that for some $n \in \N$ the objects 
$x_i\in\Z^d$, $s_i, a_i, b_i\in\mathbb{R}_+$, $\Gamma_i=[a_i,b_i)$ and the atom $A_i = (x_i,\Gamma_i)$ have already been defined for all $i \in \{0, \dots, n-1\}$. Define the set $\Delta_n := \bigcup_{i=0}^{n-1} \Gamma_i$ and remark that $\Delta_1=[0,\beta(o, 0))$. 

Given $s \in \Delta_n$, we say that the time $s$ is $n$-\emph{auspicious} if, there exist $y = y(s) \in \Z^d$ and $m = m(s) \in \{0, \dots, n-1\}$ such that
\begin{equation}\label{eq: nbom}
    s \in \I_{\langle x_m, y \rangle} \cap \Gamma_m \quad \text{and} \quad |x_m - y| \leq r_{x_m}(s);
\end{equation}
remark that in the definition above, the vertex $y = y(s) \in \Z^d$ is unique almost surely, if indeed it exists.
For each $x \in \{ x_0, \dots, x_{n-1}\}$, let
\begin{equation}
    \gamma_n(x) = \sup\{b_i : x_i=x, 0 \leq i < n\}.
\end{equation}
We say that an $n$-auspicious time $s$ is $n$-\emph{good} if one of the three situations below occurs
\begin{align}\label{eq: Situations}
    \text{S1}&) \quad  y \notin \{ x_0, \dots, x_{n-1}\}; \nonumber \\
    \text{S2}&) \quad  s \geq \gamma_n(y); \nonumber\\
    \text{S3}&) \quad s < \gamma_n(y)  \text{ ~and~ } [s,\gamma_n(y))\cap \R_{y} = \emptyset.
\end{align}

Define the set $D_n := \{s \in \Delta_n : s \text{ is $n-$good} \}$, if $D_n = \emptyset$, we stop the creation of new atoms. Otherwise, we set $s_n = \inf D_n$ and $x_n = y(s_n)$. According to \eqref{eq: Situations}, we have three possible situations: if either S1 or S2 occurs, we define $a_n = s_n$; if S3 occurs, we define $a_n = \gamma_n(y(s_n))$. To conclude our induction step, in all cases we define
\begin{equation}\label{eq: bn}
    b_n = \beta(x_n, a_n), \quad \Gamma_n=[a_n, b_n) \quad \text{and} \quad A_n = (x_n, \Gamma_n).
\end{equation}

At this point, it is interesting to note that if $s_n < a_n$, then $[s_n, a_n) \subset \Gamma_j$ for some $j \in \{0, \dots, n-1\}$ such that $x_j = x_n$.  Indeed, let $j$ be such that $b_j = \gamma_n(x_n)=a_n$, as there is at least one recovering time in the interval $(a_j,b_j)$, by condition S3, we have that $a_j<s_n$, then $[s_n, a_n) \subset [a_j,b_j) = \Gamma_j$.

Once the structure of atoms is constructed , we will partition it in generations. For each $n \geq 1$ such that $A_n$ is an atom, we say that the atom $A_{m(s_n)}$ \emph{actives} $A_n$. In this case, we denote $A_{m(s_n)} \leadsto A_n$. We define the 0-th generation as the set $\Upsilon_0 = \{A_0\}$ and for all $k \geq 1$, the $k$-th generation, $\Upsilon_k$, is the set of atoms $A$ such that 
\begin{equation}
    \exists B \in \Upsilon_{k-1} : B \leadsto A \quad \text{and} \quad  A \notin \bigcup_{i=0}^{k-1} \Upsilon_i  .
\end{equation}

For each $n \geq 0$ such that $A_n=(x_n,\Gamma_n)$ is an atom, the number of atoms activated by $A_n$, $\left|\{B : A_n \leadsto B \}\right|$, can be bounded by the random variable $X_n$ given by
\begin{equation} \label{branch}
     X_n  := \sum_{y \in \Z^d} \left| \{ t \in \I_{\langle x_n, y \rangle} \cap \Gamma_n : |x_n-y| \leq r_{x_n}(t) \}\right|.
\end{equation}

Due to the lack of memory property of the exponential distribution, the sequence of time durations for each atom $A_n$, $(b_n - a_n)_n$, is an  i.i.d. sequence with common distribution  $\beta(o,0)$. Moreover, by definitions in \eqref{eq: alphabeta} and \eqref{eq: bn}, each atom ends with a time of change, then their range $r_{x_n}(t)$, $t \in \Gamma_n$ is also independent of each other. Thus, the random variables $(X_n)_n$ defined in \eqref{branch} are i.i.d.

Hence, the stochastic sequence $(|\Upsilon_n|)_{n\geq 0}$ is dominated by a standard branching process whose number of offspring has the same distribution as $X := X_0$; that is, if this branching process dies out, then the stochastic sequence $(|\Upsilon_n|)_{n\geq 0}$ dies out too. 

The next step is to show that there exists $\lambda_0>0$ small enough such that $E[X] < 1$, that is, the branching process is subcritical. The random variable $X$ is determined only by the Poisson processes $\R_o, \T_{o}, \I_{\langle o, x\rangle}, x \in \Z^d \setminus \{o\}$ and by the random variables $N_{o,n}$, $n \geq 1$. Using the independence among these variables, we can calculate
\begin{equation}
    E[X] =  E\left[ \left| \{y \in \Z^d : |y| \leq N \} \right|\right] E[\beta(o,0)] E[Y] \leq E[(2N + 1)^d]2\lambda,
\end{equation}
where $Y\sim \exp(\lambda)$ and remark that $\beta(o,0)\sim \mbox{Gamma} (2,1)$. Since $E N^d <\infty$, we can take $\lambda_0 = 1/2E[(2N + 1)^d]>0$, then it holds that $P( |\Upsilon_n| > 0, ~\forall n \geq 0) = 0$,  $ \forall~ 0<\lambda<\lambda_0$. 

We conclude this proof showing that 
   $ \{ \zeta_t \neq \emptyset, \forall t \geq 0\} \subset \{| \Upsilon_n | > 0, \forall n \geq 0\}.$
Indeed, we will show that the collection of atoms is a  covering for all infection paths, that is, for each $x \in \Z^d$, we have $\{t > 0 : \zeta_t(x) = 1\} \subset J_x$, where $J_x = \cup_{n : x_n = x}\Gamma_n$.

Suppose that for some $m \in \N$, the sequence of vertices $(o= v_0, \dots, v_m)$ and the sequence of instants $(0 = t_0 < \dots < t_m)$ are such that, for each $k \in \{0, \dots, m-1\}$, it holds that
\begin{equation}
    t_{k+1} \in \I_{\langle v_k, v_{k+1}\rangle}, \quad 
    |v_k - v_{k+1}| \leq r_{v_k}(t_{k+1}) \quad \text{and} \quad
    [t_k, t_{k+1}] \cap \R_{v_k} = \emptyset.
\end{equation}
Arguing by induction we will show that, for each $k \in \{0, \dots, m-1\}$, the interval $[t_k, t_{k+1}]$ is covered by atoms at the vertex $v_k$, that is, $[t_k, t_{k+1}] \subset J_{v_k}$. For $k = 0$, we have that $[t_0, t_1 ] \subset \Gamma_0 \subset J_o$, since there is no recovering time in $[t_0, t_1 ]$. Suppose that for some $k  \in \{0, \dots, m-1\}$ we have that $[t_{k-1}, t_{k}] \subset J_{v_{k-1}}$. Thus, there exists $j \geq 0$ such that $t_k \in \Gamma_j$ and $x_j = v_{k-1}$. Let $\ell = \max\{ n \geq 1 : s_{n-1} \leq t_k\}$. It is clear that, $t_k \geq a_j \geq s_j$, therefore $j < \ell$. In that way, we have that $t_{k} \in \I_{\langle x_j, v_{k}\rangle} \cap \Gamma_j$ and $|v_k - x_j| \leq r_{x_j}(t_{k})$ with $j < \ell$. Thus, according to \eqref{eq: nbom}, $t_k$ is $\ell-$auspicious.

If $t_k$ is $\ell-$good, we have a contradiction: we recall that $s_{\ell} = \inf\{ s \in \Delta_{\ell} : s \text{ is $\ell-$good} \}$, thus $t_k \geq s_{\ell}$, but $\ell = \max\{ n \geq 1 : s_{n-1} \leq t_k\}$ implies that $t_k < s_{\ell}$. Thus, according to situations described in \eqref{eq: Situations}, S1, S2 and S3 do not occur, then $t_{k+1} < \gamma_{\ell}(v_k) = b_i$  for some $i < \ell$ with $x_i = v_k$. Since $i < \ell$, we have $s_i \leq t_k$, then $[t_k, t_{k+1}] \subset [s_i, b_i)$. If $s_i = a_i$ then $[t_k, t_{k+1}] \subset \Gamma_i$; if $s_i < a_i$, we write $[t_k, t_{k+1}] \subset [s_i, a_i) \cup [a_i, b_i)$ and the covering follows since $ [s_i, a_i)\subset \Gamma_{i^*}$ for some $i^* \in \{0, \dots, i-1\}$ such that $x_{i^*} = x_i$, as observed in paragraph subsequent to \eqref{eq: bn}. 
\qed

\subsection{Proof of Theorem \ref {cinf2}}\label{sub22}

\noindent 
{\textbf{Definitions and events.}}
Initially, we consider the case $d=1$. Define the sets $V=[0,6L)\cap\Z_+$, $\Delta = [0,2H)$ and let $\mathcal{B}:= V\times\Delta$ be a block in $\Z_+\times\mathbb{R}_+$, where $L$ and $H$ are integers numbers to be chosen later (as  functions of $\lambda$ and $N$). For each $v=(v_1,v_2)\in \Z_+^2$, define $\mathcal{B}_{v}:=V_{v}\times\Delta_{v}$, where $V_{v}=6Lv_1+V$ and $\Delta_{v}= 2Hv_2 + v_1+\Delta$. 

The sequence of disjoint blocks $(\mathcal{B}_{v})_{v\in\Z_+^2}$ on $\Z_+\times\mathbb{R}_+$ will induce an independent, oriented, anisotropic percolation model on a {\em renormalized graph} 
isomorphic to $\LL_+^2$, the first quadrant of the square lattice. This induced percolation processes will be defined in such a way that if there is percolation on the renormalized graph then the infection will survive in the CPDR, that is 
   $ P( \zeta_t \neq \emptyset, ~\forall t \geq 0) > 0.$

For each $v=(v_1, v_2) \in \Z^2_+$ and $k \in [0,H)\cap \Z$, define the subsets of $\Delta_v$ 
\begin{equation}\label{eq: I}
    I_v^k := 2Hv_2 + v_1 +[2k,2k + 1) \quad \text{and} \quad J_v^k:= I_v^k + 1,
\end{equation}
 and given $(j,k) \in \left( [0,6L) \times [0,H) \right) \cap \Z^2$ define the event
\begin{equation}\label{eq: Atil}
    \tilde{A}_v(j,k) := \left\{ I_v^k \cap \T_{v^j} \neq \emptyset \right\} \cap \left\{ J_v^k \cap \T_{v^j} = \emptyset \right\},
\end{equation}
 where $v^j=6Lv_1+j$. Given that $\tilde{A}_v(j,k)$ occurs, define the random variable 
   $  m := \max\{n  \geq 0 : S_{v^j,n} \in I_v^k\}, $
 where $S_{x,n}$ is defined in \eqref{eq: S}, and the event 
\begin{equation}\label{eq: A}
    A_v(j,k) := \tilde{A}_v(j,k)  \cap \left\{N_{v^j,m} \geq 7L \right\} \cap \left\{ \left(I_v^k  \cup J_v^k \right)\cap \R_{v^j} = \emptyset \right\}.
\end{equation}

By definition in \eqref{eq: R}, the event $A_v(j,k)$ implies that the range of infection of vertex $v^j$ becomes at least $7L$ at some time in $I_v^k$ and it is kept during the whole interval $J_v^k$. That is,
\begin{equation}\label{eq: P1}
    r_{v^j}(t) \geq 7L, ~\forall t \in J^k_v.
\end{equation}

The last event in the rhs of \eqref{eq: A} says that if $v^j$ is infected during the interval $I_v^k$, then the infection persists during the whole interval $J_v^k$. That is,
\begin{equation}\label{eq: P2}
 A_v(j,k) \cap   \left\{ \exists t \in I^k_v : \zeta_{t}(v^j) = 1 \right\} \subset\left\{ \zeta_{t}(v^j) = 1, ~\forall t \in J_v^k \right\}.
\end{equation}

Due to lack of memory of the exponential distribution, all the events $\{A_v(j,k) : v \in \Z^2 ,(j,k) \in \big[[0,6L)\times[0,H)\big]\cap \Z^2 \}$ are independent of each other. The occurrence of the event $A_v(j,k)$ for some coordinates $(j,k)$ will be crucial for the propagation of the infection from the block $\B_v$ to its neighbors on the renormalized graph.

For each $i \in \N$, the event $\tilde{B}_{v}^{i}(j,k)$ 
is defined as: 
\begin{equation}\label{eq: Bitil}
    \tilde{B}_{v}^{ i}(j,k) :=  \left\{ \left(I_v^k - 1\right) \cap \T_{v^j+i} \neq \emptyset \right\} \cap \left\{ I_v^k \cap \T_{v^j + i} = \emptyset \right\}.
\end{equation}
If the event $\tilde{B}^{i}_v(j,k)$ occurs, we define the random variable $m^{\ast}_i := \max\{n  \geq 0 : S_{v^j+i,n} \in I_v^k-1\}$ and the event
 \begin{equation}\label{eq: Bi}
    B^{i}_v(j,k) := \tilde{B}^{i}_v(j,k)  \cap \left\{N_{v^j+i,m_i^{\ast}} \geq i \right\}, 
\end{equation}

Observe that the variables $\{{B}^{i}_v(j,k) : i \geq 1\}$ are independent and on the event ${B}^{i}_v(j,k)$, we have that 
\begin{equation}\label{eq: P3}
    r_{v^j+i}(t) \geq |(v^j + i) - v^j|, ~\forall t \in I_v^k.
\end{equation}

Given a positive integer $M < L$ to be defined later, define the following event 
\begin{equation}\label{eq: B}
    B_v(j,k) := \left\{ \sum_{i=1}^L \ind_{\{B^{i}_v(j,k)\}} \geq M \right\}. 
\end{equation}
The occurrence of the event $ B_v(j,k)$ means that there are at least $M$ vertices in the right of $v^j$ that are able to try to infect $v^j$ during the interval $I_v^k$, these candidates will be crucial to propagate the infection inside $\B_v$.

Given $v=(v_1, v_2) \in \Z^2_+$, the definition of the next events will depend on the parity of $v_2$, the purpose is to avoid dependence in the percolation processes to be defined on the renormalized graph. If $v_2$ is even, for all $j \in [L, 2L) \cap \Z$ and $k \in [0,H)\cap \Z$, we define $C_v(j,k)$ as
\begin{equation}
    C_v(j,k) := A_v(j,k) \cap B_v(j,k) \cap \left\{ \sum_{m=0}^{H-1} \ind_{ \{A_v(\ell,m)\} } = 0, ~\forall \ell \in  [L,j)\cap \Z \right\}.
\end{equation}
The three events above are independent, since $A_v(j,k)$ uses only variables indexed by $v^j$, $ B_v(j,k)$ and $\left\{ \sum_{m=0}^{H-1} \ind_{ \{A_v(\ell,m)\} } = 0, ~\forall \ell \in  [L,j)\cap \Z \right\}$ use only variables indexed by vertices to the right of $v^j$ and to the left of $v^j$, respectively. If $C_v(j,k)$ occurs, we also define 
\begin{equation}\label{eq: V-}
    \mathcal{V}_v(j,k) := \{ v^j + i : B_v^i(j,k) \text{ occurs} \}\subseteq [v^j+1,v^j+L].
\end{equation}
Now, if $v_2$ is odd, for each $j \in [4L, 5L) \cap \Z$ and $k \in [0,H)\cap \Z$, the set $\mathcal{V}_v(j,k)$ is defined as above and the event $C_v(j,k)$ is analogously defined 
\begin{equation}\label{eq: C-}
    C_v(j,k) := A_v(j,k) \cap B_v(j,k) \cap \left\{ \sum_{m=0}^{H-1} \ind_{ \{A_v(\ell,m)\} } = 0, ~\forall \ell \in  [4L,j)\cap \Z \right\},
\end{equation}
as in the case $v_2$ even, the three events above are also dependent. 

It will be explained below that, if the event $C_v(j,k)$ occurs, on our arguments, the space-time region $\{v_j\} \times I_v^k$ is a key region inside $\B_v$ for the spread of the infection.

\vspace{.3cm}
\noindent\textbf{The Renormalization Procedure.} Once all events above are defined , we start to define the percolation on the renormalized graph. The renormalized graph will be the first quadrant of the square lattice $\LL^2_+=(\Z^2_+, \E_h\cup\E_v)$, where $\E_h = \{ (v, v+\vec{e}_1) : v \in \Z^2_+\}$, is the set of horizontal bonds, and $\E_v = \{ ( v, v+\vec{e}_2 ): v \in \Z^2_+ \}$ is the set of vertical bonds. On the set $\E_h\cup\E_v$, we define the partial order $\prec$ satisfying:

\begin{itemize}
\item $e\prec f,\ \forall e\in\E_h$ and $\forall f\in\E_v$;
\item Given $e=\left((u_1,u_2),(u_1+1,u_2)\right)$ and $f=\left((v_1,v_2),(v_1+1,v_2)\right)\in\E_h$, $e\prec f$ if and only if $u_2<v_2$ or ($u_2=v_2$ and $u_1<v_1$);
\item Given $e=\left((u_1,u_2),(u_1,u_2+1)\right)$ and $f=\left((v_1,v_2),(v_1,v_2+1)\right)\in\E_v$, $e\prec f$ if and only if $u_2<v_2$ or ($u_2=v_2$ and $u_1<v_1$).
\end{itemize}

Now, we are ready to define an exploration algorithm, concerning special events associated to the renormalized graph $\LL^2_+$. Inductively, we build a random sequence $(C_n,E_n)_{n \geq 0}$, where $C_n \subset \Z^2_+$ is the cluster of the origin $o\in\LL_+^2$ and $E_n \subset \E_h \cup \E_v$ is the set of checked bonds up to step $n$. We will also define the functions 
\begin{equation}
    \tau : \cup_{n\geq 0} C_n \longrightarrow [0,6L) \times [0,H),
\end{equation}
that associates for each vertex $v\in\cup_{n\geq 0} C_n $ coordinates in an appropriate block, these coordinates will be labels of the intervals out of which the infection spreads in the CPDR;
\begin{equation}
    \sigma: \cup_{n\geq 0} C_n \longrightarrow \Z_+
\end{equation} 
that indicates the infected vertices in the original lattice $\G^1$; and 
\begin{equation}
    \ta: \cup_{n\geq 0} C_n \longrightarrow [0,\infty)
\end{equation} 
will be defined in such a way that $\zeta_{\ta(v)}(\sigma(v)) = 1$.

To initialize the exploration process, let $o_* = (0,0) \in \Z^2_+$ be the origin of $\LL_+^2$. If  $A_{o_*}(0,0) 
$ occurs, define
\begin{equation}
    C_0 = \{o_*\}, \quad E_0 = \emptyset, \quad \tau(o_*)=(0,0), \quad  \sigma(o_*)=o \quad \text{and}  \quad \ta(o_*)=1,
\end{equation}
and observe that 
\begin{equation}
    A_v(\tau(v)) \text{ occurs \quad and \quad}  \zeta_{\ta(v)}(\sigma(v)) = 1, ~\forall v \in C_0. 
\end{equation}
Otherwise, the exploration does not start and set $(C_n, E_n) = (\emptyset, \emptyset)$ for all $n \geq 0$. 

As our induction hypothesis, suppose that for some $n \geq 0$, the sets $C_n, E_n$ are defined and for all $v \in C_n$ are also defined $\tau(v) = (\tau_1(v), \tau_2(v)),$ $\sigma(v), \ta(v)$, satisfying the conditions
\begin{equation}\label{eq: PInd}
    \sigma(v) = v^{\tau_1(v)}, \quad \ta(v) = \inf J^{\tau_2(v)}_v, \quad 
    A_v(\tau(v)) \text{ occurs}, \quad \zeta_{\ta(v)}(\sigma(v)) = 1,
\end{equation}
and $\tau_1(v) \in [0,2L)\cap \Z$ if $v_2$ is even, and $\tau_1(v) \in [3L,5L)\cap \Z$ if $v_2$ is odd.  

Define the set $F_{n} = \{ ( v,u )\in\E_h\cup\E_v : v \in C_n \text{ and } u \notin C_n \} \cap E_n^c$. If $F_n = \emptyset$, then we stop the algorithm and set $(C_m, E_m) = (C_n, E_n),\ \forall m > n$. Otherwise, let $( v, u ) \in F_n$ be the minimal bond following the order $\prec$, with $v = (v_1, v_2) \in C_n$ and  $u = (u_1, u_2) \notin C_n$ and define $E_{n+1} = E_n \cup \{ ( v,u ) \}$. To define $C_{n+1}$ we will consider four different cases. 

The first case is when $( v,u ) \in\E_h$ and $v_2$ is even. As $v$ is fixed we write $\tau(v) = (\tau_1, \tau_2)$. 
We declare that the bond $( v, u )$ is {\em open } in $\LL^2_+$ if, for some $j \in [0,L) \cap \Z$ the following event occurs:
\begin{equation}\label{eq: aberto1}
    A_u(j,\tau_2) \cap \left\{ \I_{\langle \sigma(v), u^j \rangle} \cap J_v^{\tau_2} \neq \emptyset \right\}. 
\end{equation}    

In this case, we set 
\begin{equation}
      \tau(u) = (j,\tau_2), \quad \sigma(u) = u^j \quad \text{and} \quad \ta(u) = \sup J_v^{\tau_2}.
\end{equation}
Observe that if $( v,u )$ is open, it holds $v_2 = u_2$ and $\tau_2(u) = \tau_2(v) = \tau_2$, keeping the same parity. Note also that $u_1 = v_1 + 1$, thus, by definitions in \eqref{eq: I}, we have that $ \ta(v)+ 1 = \ta(u) =  \sup J_v^{\tau_2} = \inf J_u^{\tau_2}$. By definition, $A_u(\tau(u)) = A_u(j,\tau_2)$ occurs. By the induction hypothesis, we have that $A_v(\tau(v))$ occurs and $\zeta_{\ta(v)}(\sigma(v)) = 1$, thus \eqref{eq: P1} and \eqref{eq: P2} imply that $r_{\sigma(v)}(t) \geq 7L > |\sigma(v) - \sigma(u)|,\ \forall t \in J_v^k$ and $\zeta_{t}(\sigma(v)) = 1$. By \eqref{eq: P2} and observing that the event in \eqref{eq: aberto1} occurs, we conclude that $\zeta_{\ta(u)}(\sigma(u)) = 1$. Summarizing, all properties assumed in \eqref{eq: PInd} are also true for $u$.

The case $( v,u ) \in \E_h$ and $v_2$ odd is treated analogously, we only need to replace $j \in [0,L) \cap \Z$ by $j \in [3L,4L) \cap \Z$. Clearly, all properties above are kept.

Now, suppose that $( v,u ) \in \E_v$ and $v_2$ is even. According to the order $\prec$, it is possible that for some time interval $\Delta_u$, the block $\B_u$ has already been analyzed in the CPDR. This occurs if $w := (u_1-1,u_2)\in C_n$, it means that $( w,u )$ is not open and belongs to $E_n$. More precisely, recalling that $u_2 = v_2 + 1$ is odd, the events $A_u(j,\tau_2(w))$, $j = 3L, \dots, 4L-1$ have already been analyzed. As we will see below, to define what the meaning of the event $( v,u )$ is open, we need to analyze only the events concerning the vertices $u^j$, $j=4L, \dots, 6L-1$; this will ensure that the status open or not among the analyzed bonds are independent. 

We say that the bond $( v,u )$ is open if, for some $j \in [4L, 5L) \cap \Z$ and $k \in  [0,H) 
 \cap \Z $, the event $C_u(j,k)$ occurs and for some $x \in \mathcal{V}_u(j,k)$ the following event also occurs 
\begin{equation}\label{eq: aberto2}
D_x:=    \Big\{\left( \Delta_v \cup \Delta_u \right) \cap \R_x = \emptyset \Big\}  \cap \left\{\I_{\langle \sigma(v), x \rangle} \cap J_v^{\tau_2} \neq \emptyset \right\} \cap \left\{\I_{\langle x, u^j \rangle} \cap I_u^{k} \neq \emptyset \right\}.
\end{equation}

In this case, we set 
\begin{equation}
      \tau(u) = (j,k), \quad \sigma(u) = u^j \quad \text{and} \quad \ta(u) = \inf J_v^{k}.
\end{equation}
Remark that $u_2 = v_2 +1$ is odd and $\tau_1(u) = j \in [3L, 5L)\cap \Z$. By definition, all properties in \eqref{eq: PInd}, except $\zeta_{\ta(u)}(\sigma(u)) = 1$, had been already settled. Take some $s \in \I_{\langle \sigma(v), x \rangle} \cap J_v^{\tau_2}$, the hypothesis assumed by $v \in C_n$ imply that $\zeta_s(\sigma(v)) = 1$ and $r_{\sigma(v)}(s) \geq 7L > |\sigma(v) - x|$, then $\zeta_s(x) = 1$. The first condition in \eqref{eq: aberto2} implies that $\zeta_t(x) = 1,$ $\forall t \in \left( \Delta_v \cup \Delta_u \right)\cap[s,\infty)$, in particular $\zeta_t(x) = 1,\ \forall t \in I_u^k$. Thus, under the occurrence of $C_u(j,k)$ and the third event in \eqref{eq: aberto2} (remember the definitions given in \eqref{eq: C-} and \eqref{eq: V-}, and the property \eqref{eq: P3}), we have that $r_x(t) \geq |x-u^j|,\ \forall t \in I_u^k$, then the fact $\zeta_{\ta(u)}(\sigma(u)) = \zeta_{\ta(u)}(u^j) =1$ follows from  \eqref{eq: P2}.
 
The case $( v,u ) \in \E_v$ and $v_2$  odd is analogous, replacing $j \in [4L, 5L) \cap \Z$ by $j \in [L,2L) \cap \Z$.

We emphasize that both events $( v - \vec{e}_2, v )$ is open and $( v, v+ \vec{e}_2)$ is open, accordingly first event in rhs of \eqref{eq: aberto2}, use variables indexed by vertices in $\B_v$. However, there is no dependence, since $\mathcal{V}_v(j,k) \subseteq [v^j+1,v^j+L]$ and due the parity of $v_2$, the pairs of sets like $\mathcal{V}_v(\cdot,\cdot)$ and $\mathcal{V}_{v+e_2}(\cdot,\cdot)$ are located in different halves of $\B_v$.

To conclude our inductive construction, define
\begin{equation}\label{eq: Cnum}
   C_{n+1} =
   \begin{cases}  C_n \cup \{u\}, &\text{if } ( v,u ) \text{ is open};\\
                  C_n, & \text{otherwise}.
    \end{cases}
\end{equation}

Denote $C = \cup_{n\geq 0} C_n$. As remarked above, whenever $( v,u )$ is open, it follows that the properties in  \eqref{eq: PInd} are also satisfied by $u$. In particular, our induction procedure shows that
\begin{equation}
   \zeta_{\sigma(v)}(\ta(v)) = 1, ~ \forall v \in C.
\end{equation}
Thus, if $\left\vert C\right\vert = \infty$, for all $a \in \mathbb{R}$ there exists $v = (v_1,v_2) \in C $ such that $v_1 + v_2 \geq a$, as $\ta(v) > v_1 + v_2 $ and $\zeta_t = \emptyset$ is a absorbing state, it holds that $\zeta_{\sigma(v)}(\ta(v)) = 1$ implies that $\{ \zeta_t \neq \emptyset, ~\forall t \in [0,a]\}$. As $a \in \mathbb{R}$ was chosen arbitrarily, we have that
\begin{equation}
   \left\{ \vert C \vert = \infty \right\} \subset\{ \zeta_t \neq \emptyset, ~\forall t \geq 0\}.
\end{equation}
Thus, to conclude this proof it is enough to show that the constants $L$ and $H$ can be chosen in such a way that $P \left( \vert C \vert = \infty \right) > 0$. 

\vspace{.2cm}
\noindent\textbf{Probability bounds and conclusion.} Given $v \in \Z^2_+$ and $(j,k) \in \left( [0,6L) \times [0,H) \right) \cap \Z^2$, by equations \eqref{eq: Atil} and \eqref{eq: A}, it holds that
\begin{equation}\label{eq: ProbA}
    \rho_L:=P(A_v(j,k)) = (1 - e^{-1})  e^{-1}   P(N \geq 7L)  e^{-2}.
\end{equation}
And for each $i \in \N$ and each $j \in \big[[L,2L) \cup [4L,5L)\big]\cap \Z$, the equations \eqref{eq: Bitil} and \eqref{eq: Bi} ensure that
\begin{equation}\label{eq: ProbBi}
  P(B_v^{i}(j,k))  =  (1-e^{-1})e^{-1}P(N \geq i).
\end{equation}
Thus, the hypothesis $\limsup_{n \to \infty} nP(N \geq n) > 0$ implies that 
\begin{equation}
    \sum_{i \geq 1} P(B_v^{i}(j,k)) = \infty.
\end{equation}  

From now on, fix an an arbitrary $\delta> 0$. Recall the definition given in \eqref{eq: B}. The Borel-Cantelli's Lemma implies that given any $M \in \N$, there exists $L_0(\delta, M)$ large enough such that
\begin{equation}\label{eq: ProbB}
    P(B_v(j,k)) > 1-\delta, ~\forall L \geq L_0.
\end{equation}

To analyze the probabilities of events related to vertical connections, we need to make distinction on the parity of second coordinate. Let us consider $v=(v_1, v_2) \in \Z^2$ with $v_2$ odd (the other case is analogous) and let $A_v$ be the union of independent events
\begin{equation}\label{eq: Av}
    A_v:=\bigcup_{j=4L}^{5L-1}\bigcup_{k=0}^{H-1} A(j,k).
\end{equation}
Under the occurrence of $A_v$, define 
\begin{equation}\label{eq: Jtil}
    \tilde{j} := \max \left\{ j ~:~ \sum_{k=0
    }^{H-1} \ind_{ \{A_v(\ell,k)\} } = 0, ~\forall \ell \in  [4L,j)\cap \Z \right\}
\end{equation}
and $\tilde{k} := \inf\{k :  A_v(\tilde{j}, {k}) \text{ occurs}\}$. By \eqref{eq: C-},  \eqref{eq: ProbA}, \eqref{eq: Jtil} and the bound \eqref{eq: ProbB}, we conclude that
\begin{align}\label{eq: ProbC}
    P&\left( \bigcup_{j=4L}^{5L-1}\bigcup_{k=0
    }^{H-1} C_v({j,k}) \right) \nonumber\\ &\geq \sum_{j=4L}^{5L-1}\sum_{k=0}^{H-1} \left[ P\left(A_v \cap \{(\tilde{j}, \tilde{k}) = (j,k)\} \right) P\left(C_v(j,k) \big\vert A_v \cap \{(\tilde{j}, \tilde{k}) = (j,k)\} \right) \right] \nonumber\\
    &= \sum_{j=4L}^{5L-1}\sum_{k=0}^{H-1} \left[  P\left(A_v \cap \{(\tilde{j}, \tilde{k}) = (j,k)\} \right) P\left(B_v(j,k) \right) \right] \nonumber\\
    &>(1-\delta)  P(A_v)
    = (1-\delta) \left[ 1 - (1-\rho_L)^{LH} \right], ~\forall L \geq L_0.
\end{align}

Recall the definition of $\rho_L$ in \eqref{eq: ProbA}, by the hypothesis $\limsup_{n \to \infty} nP(N \geq n) > 0$, there exists $a = a(N) > 0$ such that $\limsup_{L \to \infty} L\rho_L > a$. Fix $H=H(a,\delta)$ large enough, such that if $\rho_L > a/L$ then
\begin{equation}\label{eq: ProbC2}
\left[ 1 - (1-\rho_L)^{LH} \right] > (1-\delta).
\end{equation}

We fix $M=M(H,\lambda, \delta)$ large enough, such that under the occurrence of the event $C_v(j,k)$, we have that
\begin{equation}\label{eq: ProbM}
 P\left(\bigcup_{x \in \mathcal{V}_u(j,k)}D_x \right)\geq  1 - \left[1 - \left( e^{-4H} (1-e^{-\lambda})^2 \right) \right]^M > 1-\delta,
\end{equation}
where $D_x$ is the event defined in \eqref{eq: aberto2}.

Remark that in each step of the algorithm, all explored edges are open independently of each other.

Since for some $j \in [0,L)$ or some $j \in [3L,4L)$ (according to the appropriate parity) the event in \eqref{eq: aberto1} occurs, by equation \eqref{eq: ProbA}, if $e$ is an explored horizontal bond, that is $e \in \E_h \cap \left( \cup_{n \geq 0} E_n\right)$, it holds that 
\begin{equation}
    p_{h,L}:=  P(e~\text{is open}) = 1 - (1 - (1-e^{-\lambda})\rho_L)^{L}.
\end{equation}
Therefore, we can choose $\alpha = \alpha(\lambda, a)$ such that if $\rho_L > a/L$, it holds that $p_{h,L} > \alpha$.

Analogously, let $e \in \E_v \cap \left( \cup_{n \geq 0} E_n\right)$ be an explored vertical bond. According to conditions given in the paragraph of \eqref{eq: aberto2}, by inequalities in  \eqref{eq: ProbC} and \eqref{eq: ProbM} it holds that 
if, $\rho_L > a/L$ and $L \geq L_0$ then
\begin{equation}\label{eq: ProbV}
    p_{v,L}:= P(e~\text{is open}) > \left[ 1 - (1-\rho_L)^{LH} \right](1-\delta)^2.
\end{equation}
Combining this last inequality with \eqref{eq: ProbC2}, it holds that
\begin{equation}
    p_{v,L} > (1-\delta)^3.
\end{equation}
 
Thus, for $L$ such that $\rho_L > a/L$ and $L \geq L_0$, the exploration processes dominates  an oriented, anisotropic and independent percolation model on $\LL^2_+$,  where each horizontal or vertical bond is open with probability $\alpha$ and $(1-\delta)^3$, respectively. By a Peierls type argument we can choose $\delta=\delta(\alpha)>0$ small enough such that this anisotropic percolation model is supercritical (see Application (ii) of Theorem 3.1 in~\cite{Ke} and Theorem 3.17 in~\cite{Gr}, for an use of a Peierls type argument for anisotropic and oriented percolation, respectively ). This concludes the proof for the case $d=1$.

For $d \geq 2$ the proof is essentially the same. The main difference is on the family of blocks $(\B_v)_{v \in \Z_+^2}$. We set $V=\left[ [0,6L)\times[0,L)^{d-1} \right] \cap \Z^{d}$, $\Delta = [0,2H)$ and let $\mathcal{B}= V\times\Delta$ be a block in $\Z_+^{d} \times \mathbb{R}_+$. For each $v=(v_1,v_2)\in \Z_+^2$, define $\mathcal{B}_{v}=V_{v}\times\Delta_{v}$, where $V_{v}=6Lv_1\vec{e}_1+V$ and $\Delta_{v}= 2Hv_2 + v_1+\Delta$.  Analogously, throughout the family $(\B_v)_{v \in \Z_+^2}$, we construct a coupling between the CPDR on $\G^d$ and a supercritical percolation model on $\mathbb{L}^2_+$.  The main difference on the probability bounds of vertical and horizontal connections occurs on the expression $[1 - (1- \rho_L)^{LH}]$ in \eqref{eq: ProbC}, \eqref{eq: ProbC2} and \eqref{eq: ProbV}. In the general context, this expression turns into $[1 - (1- \rho_L)^{L^d H}]$ and the hypothesis $\limsup_{n \to \infty} nP(N^d \geq n) > 0$ yields the key fact $\limsup_{L \to \infty} L^{d}\rho_L > 0$.  

Remark that the CPDR survive with positive probability on a \emph{slab} of $\G^d$, that is, on the subgraph induced by the set of vertices $\left[ \Z \times [0,L)^{d-1} \right] \cap \Z^d$. 
\qed

\section{Anisotropic percolation with random range}\label{s3}

\noindent The proof of Theorem~\ref{perc1} has some common features with the proof of Theorem~\ref{cfinito2} like the domination by a subcritical branching process, but the details are different. The proof of Theorem~\ref {perc2}, although much more involved than that Theorem~\ref {perc1}, follows the same script as the proof of Theorem~\ref {cinf2} (indeed the former is simpler than the latter). 

We remark that the statement of Theorem~\ref{perc2} is easily proved for non-oriented percolation. Indeed, under the hypothesis $\sum_{n\in\N} P(N\geq n)= +\infty$, by Borel-Cantelli's Lemma, for each vertex $v$ and for all $p\in (0,1]$, the set 
$\{n \in \N : \langle -n\vec{e}_1,o\rangle \text{ and } \langle -n \vec{e}_1, v \rangle \text{ are open}\}$
is infinite $\emph{a.s.}$

\subsection{Proof of Theorem~\ref {perc1}}\label{sub31}

\noindent We will consider $d=2$, the same proof can be carried out with very minor modifications for all values of $d > 2$, indeed only modifications on the notation are needed. Given $l\in\N$ consider the partition $(\mathcal{S}^l_v)_{v\in\Z^2}$, of $\Z^2$ in one-dimensional boxes $\mathcal{S}^l_{(v_1,v_2)}:= (lv_1,v_2)+ \big[ \{0,1,\dots,l-1\}\times \{0\} \big]$. This partition will induce a percolation process in the renormalized graph $\tilde{G}=(\Z^2,\E_v\cup \E_h)$ isomorphic to $G$. We remark that for $d>2$ we still consider a partition of $\Z^d$ in one-dimensional boxes like $\mathcal{S}^l_{(v_1, \dots, v_d)}:= (lv_1,v_2, \dots, v_d)+ \big( \{0,1,\dots,l-1\}\times \{0\}^{d-1} \big)$.

We say that a bond $(u,v)$ in the renormalized lattice $\tilde{G}$ (that is $v$ is $u+\vec{e}_2$ or $u+ n\vec{e}_1$ for some $n\in\N$) is {\em good} if and only if the event occurs 
\begin{equation}
    \bigcup_{x\in \mathcal{S}^l_u}\bigcup_{y\in \mathcal{S}^l_v}\left\{(x,y)\mbox{ is open in } G\right\},
\end{equation} 
this construction ensures that 
\begin{equation}\label{renorm}
\{(0,0)\rightarrow\infty\mbox{ in }G\}\subset\{(0,0)\leadsto\infty\mbox{ in }\tilde{G}\},
\end{equation}
here we are using the notation $\leadsto$ to denote connections by good paths in the renormalized graph. 

If $(u,v)$ is a vertical bond of $\tilde{G}$, that is $(u,v)\in\E_v$, then 
\begin{equation}\label{cotav}
P((u,v)\mbox{ is good})=1-(1-q)^l.
\end{equation} 

Given $p < 1$, if $(u,v)\in\E_h$ with $v=u+\vec{e}_1$, we claim that there exists some constant $c(p)<1$, not depending on $l$ or $q$, such that 
\begin{equation}\label{cotah1}
P((u,v)\mbox{ is good})<c(p)
\end{equation}
and from now on we fix $l$ such that 
\begin{equation}\label{defl}
\sum_{n\geq l}P(N\geq n)<\frac{1-c(p)}{2}.
\end{equation}
We postpone the demonstration of this claim to the end of this proof. If $(u,v)\in\E_h$ with $v=u+n\vec{e}_1$,  $n\geq 2$, it holds that 
\begin{equation}\label{cotahn} 
P((u,v)\mbox{ is good}) 
\leq\sum_{i=1}^l P(N_{(lu_1+l-i,u_2)}\geq l(n-1)+i)
=\sum_{i=1}^l P(N\geq l(n-1)+i).
\end{equation}

We can use the bounds (\ref{cotav}), (\ref{cotah1}) and (\ref{cotahn}) to estimate the expected number of neighbours of any vertex $v$ of $\tilde{G}$ in its good cluster:
\begin{align}\label{cotaviz}
\nonumber E |\{u\in\Z^2 : (v,u) \mbox{ is  good} \}| 
&< 1-(1-q)^l + c(p) +\sum_{n\geq 2}\sum_{i=1}^l P(N\geq l(n-1)+i)\\
&\leq 1-(1-q)^l + c(p) +\sum_{i\geq l} P(N\geq i).
\end{align}

Due to inequalities (\ref{defl}) and (\ref{cotaviz}), we can take $q$ close enough to 0 such that
$E |\{u\in\Z^2 : (v,u) \mbox{ is good} \}| < 1.$
Thus, the percolation processes on the renormalized graph $\tilde{G}$ is dominated by a subcritical branching processes, then by inclusion (\ref{renorm}), we have that $\theta(p,q)=0$, yielding that $q_c(p)>0$, $\forall p<1$.

We conclude proving the claim made above: Given any $p<1$ and $l\in\N$ there exists a constant $c:=c(p)<1$ not depending on $l$ (indeed $c$ depends on the distribution $N$), such that
\begin{equation}
    P((u,v)\mbox{ is good}) < c(p),\ \forall (u,v)\mbox{ with } v = u+\vec{e}_1.
\end{equation}
For this, by translation invariant, it is enough to prove that $P(A)<1$, where $A$ is defined as follows:
\begin{equation}
    A:=\bigcup_{i=1}^\infty\bigcup_{j = 0}^{\infty}\Big((-i\vec{e}_1, j\vec{e}_1)\mbox{ is open in }G\Big).
\end{equation}
Define $i_0:=\inf\{i \geq 0 : P(N = i) > 0\}$. Let $A_0 = \emptyset$ and, for $i_0 \geq 1$, consider the event $A_{i_0}$ defined by 
\begin{equation}
    A_{i_0}:=\bigcup_{i=1}^{i_0}\bigcup_{j =  0}^{\infty}\Big((-i\vec{e}_1,j \vec{e}_1 )\mbox{ is open in }G\Big).
\end{equation} 
Observe that, since $p < 1$, it holds that $P(A_{i_0}) < 1$. Therefore  
\begin{equation}
P(A^c) \geq P\left(\bigcap_{i> i_0} \{N_{(-i,0)}<i\}\cap A_{i_0}^c\right)
= \prod_{i > i_0} \left[ 1-P(N\geq i) \right] \cdot P\left(A_{i_0}^c\right)>0,
\end{equation}
since $E N=\sum_{i=1}^\infty P\{N\geq i\}<+\infty$. This concludes the proof of our claim. 
\qed

\smallskip
\begin{remark}
Observe that if $i_0=0$, that is $P(N=0)>0$, we can write 
\begin{equation}
P(A^c) \geq P\left(\bigcap_{i> 0} \{N_{(-i,0)}<i\}\cap A_{i_0}^c\right)
= \prod_{i > 0} \left[ 1-P(N\geq i) \right] >0.
\end{equation}
Thus, as the equation above does not depend on $p$, we can extend Theorem~\ref{perc1} for $p=1$.
\end{remark}

\subsection{Proof of Theorems~\ref {perc2}}
\label{sub32}
\noindent We will prove that $P(o\rightarrow\infty)>0$ for all $p,q>0$. It is enough to show this only for $d=2$, since the APRR in $d+1$ dimensions can be embedded in the $d$-dimension case for $d\geq 2$. 

We follow the same script of the proof of Theorem~\ref{cinf2}. Since we are considering oriented percolation, we have some changes on the geometry of the renormalized lattice. Now, the family of blocks $(\mathcal{B}_{v})_{v\in\Z_+^2}$ of $\Z_+^2$ is defined as follows: given integer numbers $L$ and $H$, for each vertex $v = (v_1, v_2)$ of the lattice $\mathbb{L}_+^2$, let $o(v) = (3Lv_2 + 6Lv_1, Hv_2) \in \Z^2$ and define $\B_v = o(v) + \B$, where $\B = \big[ [0,6L) \times [0,H) \big] \cap \Z^2$.

For each $v \in \Z_+^2$ and each $(j,k) \in \B$, define $v(j,k) = o(v) + (j,k)$ the associated vertex in $\B_v$. Fix $M \leq L$ an integer number, we also define the following events
\begin{align}\label{eq: eventosnovos}
   & A_v(j,k) = \{ N_{v(j,k)} \geq 7L \};  \nonumber \\
& B^i_v(j,k) = \{ N_{v(j,k) - i \vec{e_1} } \geq i \}, ~i \geq 1\mbox{ and } \nonumber \\
& B_{v}(j,k) = \{\big\vert \mathcal{V}_v(j,k) \big\vert \geq M\}, \text{ where } \mathcal{V}_v(j,k) = \{ i \in [1,L] :  B^i_v(j,k) \text{ occurs} \}.
\end{align}
Moreover, for each $v \in \Z_+^2$ if $(j,k) \in \B$ is such that $j \in [2L, 3L)$, define
\begin{equation}
    C_v(j,k) = A_v(j,k) \cap B_v(j,k) \cap \left\{ \sum_{m=0}^{H-1} \ind_{ \{A_v(\ell,m)\} } = 0, ~\forall \ell \in  (j,3L)\cap \Z \right\}.
\end{equation}

Now, we define the exploration algorithm on the oriented lattice $\LL^2_+$. Inductively, we build a random sequence $(C_n,E_n)_{n \geq 0}$, where, we recall, $C_n \subset \Z^2_+$ and $E_n \subset \E_h \cup \E_v$ is the set of checked bonds up to step $n$. We will also define the \emph{golden coordinates}
\begin{equation}\label{eq: goldtau}
    \tau : \cup_{n\geq 0} C_n \longrightarrow [0,3L) \times [0,H),
\end{equation}
that associates for each vertex $v \in \cup_{n\geq 0} C_n$, the \emph{golden vertex} 
\begin{equation}\label{eq: goldv}
    v^{\ast} = v(\tau(v)) \in \B_v.
\end{equation} 
By construction we will have that, if $v \in \cup_{n \geq 0} C_n$ then $\big(o \to v^{\ast} \big)$ occurs on the APRR model.

If  $A_o(0,0)$ occurs, define $C_0 = \{o\}$, $E_0 = \emptyset,$ and $\tau(o)=(0,0)$. Observe that, according to \eqref{eq: eventosnovos}, \eqref{eq: goldtau} and \eqref{eq: goldv} we have
\begin{equation}
   N_{v^{\ast}} \geq 7L \text{~~and~~} o \to v^{\ast},  ~\forall v \in C_0. 
\end{equation}
Otherwise, the exploration does not start and set $(C_n, E_n) = (\emptyset, \emptyset)$ for all $n \geq 0$. 

As our induction hypothesis, suppose that for some $n \geq 0$, the sets $C_n, E_n$ are defined and, for all $v \in C_n$, we have also defined $\tau(v)$ and $v^{\ast}$ satisfying the conditions
\begin{equation}\label{eq: HI}
 N_{v^{\ast}} \geq 7L \text{~~and~~} o \to v^{\ast},  ~\forall v \in C_n. 
\end{equation}
Define the set $F_{n} = \{ ( v,u ) : v \in C_n \text{ and } u \notin C_n \} \cap E_n^c$. If $F_n = \emptyset$, then we stop the algorithm and set $(C_m, E_m) = (C_n, E_n),\ \forall m > n$. Otherwise, let $( v, u ) \in F_n$ be the minimal bond following the order $\prec$ (defined in Section~\ref{sub22}), with $v = (v_1, v_2) \in C_n$ and  $u = (u_1, u_2) \notin C_n$ and define $E_{n+1} = E_n \cup \{ ( v,u ) \}$. We denote $\tau(v) = (\tau_1, \tau_2)$. To define $C_{n+1}$ we will consider two different cases. 

Suppose that $(v,u) \in \E_h$. We say that $(v,u)$ is open if, for some $j \in [0,L) \cap \Z$, the following
event occurs
\begin{equation}\label{eq: openv}
    \left\{\big( v^{\ast} , u(j,\tau_2) \big) \text{ is open in } G_{\textbf{N}} \right\} \cap
    A_u(j,\tau_2).
\end{equation}
In this case we set $\tau(u) = (j,\tau_2)$. Note that we assume \eqref{eq: HI}, hence $\big( v^{\ast},u(j,\tau_2) \big)$ is indeed an edge in the random graph $G_{\textbf{N}}$. According to definition of $A_u(j,k)$ in \eqref{eq: eventosnovos}, combining \eqref{eq: HI} and \eqref{eq: openv}, it holds that $ o \to u^{\ast}$ and $N_{u^{\ast}} \geq 7L$.

Suppose that $(v,u) \in \E_v$. We say that $(v,u)$ is open if, for some $(j,k) \in \big[ [2L,3L)\times [0,H) \big]\cap \Z^2$ it holds that $C_u(j,k)$ occurs and also, there exists $i \in \mathcal{V}_u(j,k)$ such that $D_i$ occurs, where $D_i = \cap_{\ell = 1}^4 D_i^{\ell}$ with
\begin{align}\label{eq: openh}
      D_i^1 &:=  \left\{\big(v^{\ast} , v^{\ast} + (3L  - \tau_1 + j - i) \vec{e}_1 \big) \text{ is open in } G_{\textbf{N}} \right\}; \nonumber \\
 D_i^2 &:= \left\{\big(u(j,k) - i\vec{e}_1 ,u(j,k) \big) \text{ is open in } G_{\textbf{N}} \right\}; \nonumber \\
 D_i^3 &:=\bigcap_{m=0}^H \left\{\big(v(3L+j-i,m), v(3L+j-i,m)+ \vec{e}_2 \big) \text{ is open in } G_{\textbf{N}} \right\}\mbox{ and} \nonumber \\
 D_i^4 &:=\bigcap_{m=0}^{H-1} \left\{\big(u(j-i,m), u(j-i,m)+ \vec{e}_2 \big) \text{ is open in } G_{\textbf{N}} \right\}.
\end{align} 
In this case we set $\tau(u) = (j,k)$. First note that $A_u(j,k)$ occurs, thus $N_{u^{\ast}} \geq 7L$. Throughout this paragraph, as an auxiliary notation, let $x = v^{\ast} + (3L + j - \tau_1 - i)\vec{e}_1$  and $y = u(j,k) - i\vec{e}_1$. Assuming condition \eqref{eq: HI}, it holds that $(v^{\ast},x)$ is an edge in the random graph $G_{\textbf{N}}$. We also have that $(y, u(j,k))$ is an edge in the random graph $G_{\textbf{N}}$, since $B^i_u(j,k)$ occurs. Note that $x$ and $y$ have the same first coordinate, and since $D_i^3$ and $D_i^4$ occur, all edges in the column contain $x$ and $y$ and contained in $\B_v \cup \B_u$ are open, therefore $x \to y$. Wrapping up, we conclude that $o \to u^{\ast}$ and $N_{u^{\ast}} \geq 7L$.

To finish our inductive procedure, define
$C_{n+1}$ as in \eqref{eq: Cnum}. 

We emphasize that, by construction, for each explored edge $e \in \cup_{n \geq 0} \E_n$, the event $\{e \text{ is open}\}$ is independent of all others. 
   
By hypothesis $\limsup_{n \to \infty}nP(N \geq n) > 0$, there exists $a > 0$, such that $\mathcal{L} = \{ L \in \N : P(N \geq 7L) > a/L \}$ is infinity. Let $\alpha= \alpha(N,p)$ be such that 
\begin{equation}\label{eq: alpha2}
   1-  \Big( 1- P(N \geq 7L)p \Big)^L \geq \alpha, ~\forall L \in \mathcal{L}.
\end{equation}
Fix an arbitrary $\delta > 0$. Let $H = H(N,\delta)$ be such that
\begin{equation}
   1-  \Big( 1- P(N \geq 7L) \Big)^{LH} \geq (1-\delta), ~\forall L \in \mathcal{L}.
\end{equation}
Fix $M = M(H,p,q, \delta)$ such that
\begin{equation}
   1-  \Big( 1- p^2 q^{4H-1} \Big)^{M} \geq (1-\delta), ~\forall L \in \mathcal{L}.
\end{equation}
Since $\sum_{n \geq 1} P(N \geq n) = \infty$, there exists $L_0 = L_0(N,M,\delta)$ such that 
\begin{equation}\label{eq: Lzero}
    P(B_v(j,k)) \geq (1-\delta), ~\forall L \geq L_0.
\end{equation}

Thus, according to conditions given in paragraphs of \eqref{eq: openv} and \eqref{eq: openh}, if $L \in \mathcal{L}$ and $L \geq L_0$ the inequalities \eqref{eq: alpha2} -- \eqref{eq: Lzero} ensure that, for explored edges,  $P(e \text{ is open}) \geq \alpha$, if $e \in \E_h$, and (arguing as in the paragraph of  \eqref{eq: ProbC}) $P(e \text{ is open}) \geq (1 - \delta)^3$, if $e \in \E_v$. As $\delta$ was fixed arbitrarily, the theorem follows if we take $\delta$ sufficiently small.
\qed

\subsection{Proof of Theorem~\ref {perc4}}
\label{sub33}
\noindent We define recursively the sequence $(j_n)_{n\geq 0}$ as follows: set $j_0 = o$ and $u_0 = N_0$. Given $n \geq 0$ such that $j_n$ and $u_n$ are already defined, for each $i \in (j_{n}, u_{n}]$, define  
$B_i = \{(j_{n}, i) \text{ is open}\}$ and
$\tilde{N}_i = N_i \mathbbm{1}_{B_i}$. Let $a_{n+1} = \max\{ i + \tilde{N}_i : i \in (j_{n}, u_{n}]\}$, if $a_{n+1} = u_n$, we say that $u_n + 1$ is a cutting point and restart the recursion defining $j_{n+1} = u_n + 1$, otherwise, we define $j_{n+1} \in (j_{n}, u_{n}]$ in such a way that $j_{n+1} + \tilde{N}_{j_{n+1}} = a_{n+1}$. Observe that, in the last case, $j_{n+1}$ could be not unique, to fix a criteria, we can choose the biggest one among all the possibilities. To conclude the recursion, in both cases, we define $u_{n+1} = j_{n+1}  + N_{j_{n+1}}$.

Let $A_i$ be the event where the vertex $i$ is a cutting point. Define $\tilde{N} = XN$, where $X$ is a Bernoulli random variable with parameter $p$. Given any $\delta > 0$, we can find some large constant $C>0$, such that for all $i$ large enough, it holds
\begin{equation}\label{eq: PAI2}
    P(A_i) = \prod_{\ell = 1}^{i}P(\tilde{N} < \ell) \leq C \exp\left\{ -p\beta(1-\delta) \sum_{\ell = i_o}^i \frac{1}{\ell}\right\}.\ 
\end{equation}
Thus, if $\beta > 1/p$, we can take $\delta$ small enough such that $\sum_{i\geq 0} P(A_i) < \infty$. Therefore, by Borel-Cantelli Lemma, it holds that $P(A_i ~ \text{i.o.}) = 0$.

Note that
\begin{equation}\label{eq: existeAi2}
   \{ |C_0| < \infty\} \subset \bigcup_{i \geq 1} A_i.
\end{equation}

By the strong Markov property, if $P(\cup_{i \geq 1} A_i) = 1$ then $P(A_i ~\text{i.o.}) = 1$. Thus $P(\cup_{i \geq 1} A_i) < 1$ and the result follows from \eqref{eq: existeAi2}.
\qed


\vspace{1cm}

\noindent{\Large \bf Acknowledgments:} P.A.G. was supported by S\~ao Paulo Research Foundation (FAPESP) grants 2020/02636-3 and 2017/10555-0. B.N.B.L. was supported in part by CNPq grant 305811/2018-5, FAPEMIG (Programa Pesquisador Mineiro) and FAPERJ (Pronex E-26/010.001269/2016). Both authors thank the two anonymous referees for their valuable work.

\end{document}